\documentclass[11pt,twoside,a4paper]{article}
\usepackage{amsfonts,amssymb,amsmath,latexsym,epsfig,verbatim,enumerate}

\newcommand{\kk}{\mathbb{K}}

\newcommand{\rr}{\mathbb{R}}
\newcommand{\qq}{\mathbb{Q}}

\newcommand{\al}{\alpha}
\newcommand{\fn}[1]{\frac{#1}{n}}

\newcommand{\be}{\beta}
\newcommand{\ga}{\gamma}

\newcommand{\ep}{\epsilon}

\title{On the free energy of a directed polymer in a Brownian environment}
\author{John Moriarty and Neil O'Connell \\ \small{Department of Mathematics and BCRI, }\\ \small{University College Cork, Ireland} \\ \small{Revised 13th June 2006}}
\date{}

\begin{document}
\maketitle
\bibliographystyle{plain}

\begin{abstract}
We prove a formula conjectured in \cite{MR1865759} for the free energy density of a directed polymer in  a Brownian environment in $1+1$ dimensions.

\end{abstract}


\section{Introduction}

Let $B^{(1)},B^{(2)},\ldots$ be independent standard one-dimensional
Brownian motions.
Denote the increments of $B^{(i)}$ by $B^{(i)}_{(s,t)}=B^{(i)}_t-B^{(i)}_s$.

For $\beta \in \mathbb{R}$ set
\begin{eqnarray} \label{zed}
	Z_n(\beta)=\int_{0<s_1<\ldots<s_{n-1}<n}ds_1 \ldots ds_{n-1}
	\exp \big\{ \beta(B^{(1)}_{(0,s_1)}+\ldots+B^{(n)}_{(s_{n-1},n)})
	\big\} .
\end{eqnarray}
This is the partition function for a continuous model of a directed polymer in a Brownian environment in 
$1+1$ dimensions.  In the paper~\cite{MR1865759}, using queueing-theoretic ideas
in the context of geometric functionals of Brownian motion, certain limiting results
were obtained which led the authors to conjecture an explicit formula for the
free energy density
\begin{eqnarray} \label{fre}
	\lim_{n \to \infty} \frac{1}{n} \log Z_n(\beta),
\end{eqnarray}
namely that it should be given by, almost surely, 
\begin{eqnarray}
	f(\beta)=
	\left\{ \begin{array}{ll}-(-\Psi)^*(-\beta^2)-2\log |\beta| & :\beta \neq 0\\ 1 & : \beta=0	
		\end{array} \right.
	\label{putz}
\end{eqnarray}
where $\Psi(m) \equiv \Gamma'(m)/\Gamma(m)$ is the restriction of the digamma function to $(0,\infty)$, and
$(-\Psi)^*$ is the convex dual of the function $-\Psi$.
The aim of this paper is to give a rigorous proof of this conjecture.
The proof uses tools from large deviation theory.
As a corollary we give a new proof that $c \geq 2$, where 
\begin{eqnarray}
	\lim_{n \to \infty} \frac{1}{n}L_n(n) = c \quad \mbox{ a.s.} \label{drea}
\end{eqnarray}
and
\begin{eqnarray*}
	L_n(t)=\sup_{0 \leq s_1 \leq \ldots \leq s_{n-1}\leq t}B^{(1)}_{(0,s_1)}+\ldots+B^{(n)}_{(s_{n-1},t)}.
\end{eqnarray*}
It was proved using direct methods in \cite{MR1935124} that $c=2$, where the authors also describe how the result may be deduced from the theory of random matrices.

The directed polymer model we have discussed here, and for which we have computed
the free energy density, is a continuous version of the classical two-dimensional directed
polymer, where it is not known how to compute the free energy density (see, for
example, Derrida \cite{MR1043640}, Carmona and Hu \cite{MR1939654}). Recent work on a continuous model different to this can be found in Comets and Yoshida \cite{MR2117626}.

In the next section we recall the framework which was developed in~\cite{MR1865759} to 
extend some standard constructions from queueing theory to the context of geometric 
functionals of Brownian motion, and explain how this leads to the conjectured formula
for the free energy density.   Section \ref{flint} is devoted to the
proof of the main result. In section \ref{an} we prove that $f$ is analytic and strictly convex, and record a large deviation principle that we will use in section \ref{appl} to prove (\ref{drea}).

\noindent{\em Acknowledgement:}
Research supported in part by Science Foundation Ireland,
Grant Number SFI 04/RP1/I512.  The authors would also like to
thank the anonymous referee for helpful suggestions which have
led to a much improved presentation.

\section{Generalised Brownian queues}
\label{gbqt}
In this section we recall the framework which was developed in~\cite{MR1865759} to 
extend some standard constructions from queueing theory to the context of geometric 
functionals of Brownian motion, and explain how this leads to the conjectured formula
for the free energy density.

The \emph{generalised Brownian queue} is characterised as follows.
Let $B$ and $C$ be two independent standard Brownian motions indexed by the 
entire real line, and write 
\begin{eqnarray*}
	B_{(s,t)}=B_t-B_s, \quad C_{(s,t)}=C_t-C_s
\end{eqnarray*}
Fix $m>0$ and, for $t \in \mathbb{R}$, set 
\begin{eqnarray*}
r(t)&=&\log \int^t_{-\infty} ds \exp\{B_{(s,t)}+C_{(s,t)}-m(t-s)\}\\
f(s,t)&=&B_{(s,t)}+r(s)-r(t)\\
g(s,t)&=&C_{(s,t)}+r(s)-r(t)
\end{eqnarray*}
and define $f:\rr \to \rr$ by $f(t)=f(0,t).$

To put this in context, the `Brownian queue' is defined similarly but
with `log $\int$ exp' replaced by `sup' (see, for example, \cite{MR1865759}).
The usual $M/M/1$ queue is defined similarly to the Brownian queue but
with Brownian motions replaced by Poisson counting processes.
Thus the Brownian motions $B_t$ and $mt-C_t$ can be thought of, respectively, 
as the \emph{arrivals} and \emph{service} processes, $r$ as the \emph{queue-length process} 
and $f$ as the \emph{output}, or \emph{departure}, process. 

In~\cite{MR1865759} it is shown, using results of Matsumoto and Yor~\cite{MR1696208, MR1833628},
that the generalised Brownian queue is {\em quasi-reversible}, that is:
$f$ is a standard Brownian motion, and $\{f(s), s \leq t\}$ is independent of $r(t)$.
We can thus consider a sequence of generalised Brownian queues in tandem and expect this
`queueing network' to have nice properties (analogous to the `product-form solutions' 
of classical queueing theory).

Let $B,B^{(1)}, B^{(2)},\ldots$ be a sequence of independent standard Brownian motions, 
each indexed by $\mathbb{R}$, and let $m>0$ be a fixed constant. 
For $-\infty<s \leq t < \infty$, set 
\begin{eqnarray*}
r_1(t)&=&\log \int^t_{-\infty} ds \exp\{B_{(s,t)}+B^{(1)}_{(s,t)}-m(t-s)\}\\
f_1(s,t)&=&B_{(s,t)}+r_1(s)-r_1(t)\\
g(s,t)&=&B^{(1)}_{(s,t)}+r_1(s)-r_1(t),
\end{eqnarray*}
for each $k=2,3,\ldots$ set
\begin{eqnarray*}
r_k(t)&=&\log \int^t_{-\infty} ds \exp\{f_{k-1}(s,t)+B^{(k)}_{(s,t)}-m(t-s)\}\\
f_k(s,t)&=&f_{k-1}{(s,t)}+r_k(s)-r_k(t),
\end{eqnarray*}
and for all $k$ define $f_k:\rr \to \rr$ by $f_k(t)=f_k(0,t).$

Note that $r_1(t)$ is clearly stationary in $t$; to see that $r_1(0)<\infty$ almost surely 
simply note that, with probability one, $B_{(s,0)}+B^{(1)}_{(s,0)}+ms<ms/2$ for all $s$ 
sufficiently negative (by Strassen's law of the iterated logarithm, for example). 
In fact, $r_1(0)$ has the same law as $-\log Z_m$, where $Z_m$ is gamma-distributed 
with parameter $m$: this is Dufresne's identity~\cite{MR1825324, MR1833627}. 

We first state the quasi-reversibility property, as presented in \cite{MR1865759}.

\newtheorem{ino}{Theorem}
\begin{ino} \label{inoz}
\begin{enumerate}
\item $f_1$ and $g$ are independent standard Brownian motions indexed by $\mathbb{R}$.
\item For each $t \in \mathbb{R}$, $\{(f_1(s),g(s)), -\infty<s \leq t \}$ is independent of $\{r_1(s), s \geq t\}$.
\end{enumerate}
\end{ino}
It follows from Theorem \ref{inoz} that $r_1(0),r_2(0),\ldots$ is a sequence of i.i.d. random variables, each distributed as $-\log Z_m$. By construction, we have
\begin{eqnarray} \label{dag}
\lefteqn{\sum_{k=1}^n r_k(0)=\log \bigg[ \int_{-\infty}^0 du \exp (B_{(u,0)}+mu)}\nonumber\\
&&\times \int_{u<s_1<\ldots<s_{n-1}<0}ds_1\ldots ds_{n-1} \exp\{B^{(1)}_{(u,s_1)}+\ldots+B^{(n)}_{(s_{n-1},0)}\} \bigg]
\end{eqnarray}

Applying the strong law of large numbers, it can be deduced 
(see \cite{MR1865759} for details) that:

\newtheorem{ino2}[ino]{Theorem}
\begin{ino2} \label{ino2z}
For each $m>0$:
\begin{eqnarray*}
&&\lim_{n \to \infty} \frac{1}{n}\log \int_{-\infty<u<s_1<\ldots<s_{n-1}<0} duds_1\ldots ds_{n-1}\\
&& \quad \exp\{mu+B^{(1)}_{(u,s_1)}+\ldots+B^{(n)}_{(s_{n-1},0)}\}=-\Psi(m)
\end{eqnarray*}
almost surely, where
\begin{equation*}
\Psi(m)=E \log Z_m=\Gamma'(m)/\Gamma(m)
\end{equation*}
is the digamma function (and $\Gamma$ is the Gamma function).
\end{ino2}
Theorem \ref{ino2z} can be interpreted as follows. Let $\mathcal{B}$ denote the $\sigma$-field generated by the Brownian motions $B^{(1)},B^{(2)},\ldots,$ and let $\tau_1,\tau_2,\ldots$ be the points of a unit-rate Poisson process on $\mathbb{R}_+$, independent of $\mathcal{B}$. For $t_0,t_1,\ldots,t_n \in \rr$ define 
\begin{eqnarray*}
	E_n(t_0,t_1,\ldots,t_n)&=&B^{(1)}_{(t_0,t_1)}+\ldots+B^{(n)}_{(t_{n-1},t_n)}\\
	F_n^m(t_0,t_1,\ldots,t_{n-m+1})&=&
	\exp(B^{(n)}_{(t_0,t_1)}+\ldots+B^{(m)}_{(t_{n-m},t_{n-m+1})})\\
	F_n&=&F_n^1
\end{eqnarray*}
By Brownian scaling, Theorem \ref{ino2z} is equivalent to:

\newtheorem{ino3}[ino]{Theorem}
\begin{ino3} \label{ino3z}
For $\theta \neq 0$,
\begin{equation*}
\lim_{n \to \infty} \frac{1}{n} \log E[\exp(\theta E_n(0,\tau_1,\ldots,\tau_n))|\mathcal{B}]=
-2\log|\theta|-\Psi(1/\theta^2),
\end{equation*}
almost surely.
\end{ino3}
 
Thus, if we set 
\begin{eqnarray*}
\Lambda(\theta)=\left\{ \begin{array}{ll}
-2\log |\theta| - \Psi(1/\theta^2), & \theta \neq 0 \\  
0, & \theta =0,
\end{array} \right.
\end{eqnarray*}
we have
\begin{eqnarray*}
	\lim_{n \to \infty}\frac{1}{n} \log E[\exp(\theta E_n(0,\tau_1,\ldots,\tau_n))| \mathcal{B}]=\Lambda(\theta),
\end{eqnarray*}
almost surely. From the asymptotic expansion
\begin{eqnarray} \label{coff}
	\Psi(x) \sim \log x-\frac{1}{2x}-\sum_{k=1}^{\infty}\frac{B_{2k}}{2kx^{2k}}
\end{eqnarray}
as $x \to \infty$ (see, for example, Abramowitz and Stegun \cite{MR1225604}), we have that $\Lambda$ is finite and differentiable everywhere, with $\Lambda(0)=\Lambda'(0)=0$. It follows that the sequence $\frac{1}{n}E_n(0,\tau_1,\ldots,\tau_n)$
satisfies the following conditional large deviation principle:
\newtheorem{test}[ino]{Theorem}
\begin{test} \label{testz}
Given $\mathcal{B}$, $\frac{1}{n}E_n(0,\tau_1,\ldots,\tau_n)$ satisfies a large deviation principle with good rate function
\begin{eqnarray*}
	\Lambda^*(x)=\sup_{\theta \in \mathbb{R}}[x\theta-\Lambda(\theta)]
\end{eqnarray*}
almost surely.
\end{test}

This is a \emph{quenched large deviation principle}, associated with the conditional law of large numbers. For example, Theorem \ref{testz} implies  that given $\mathcal{B}$, $\frac{1}{n}E_n(0,\tau_1,\ldots,\tau_n) \to 0$ almost surely. Another implication is that for any $x>0$,
\begin{eqnarray*}
	\lim_{n \to \infty}\frac{1}{n}\log P(E_n(0,\tau_1,\ldots,\tau_n)>xn|\mathcal{B})=-\Lambda^*(x)
\end{eqnarray*}
almost surely.
For two other related large deviation principles see \cite{MR1865759}.

We will now describe how this relates to the Brownian directed polymer model.  
It is shown in Lemma \ref{km} using Kingman's subadditive ergodic theorem that there exists a function $\ga:\rr \to \rr$ such that given $x<0$, 
\begin{eqnarray} \label{gam}
\lim_{n \to \infty} \frac{1}{n}\log \int_{xn<s_1<\ldots<s_{n-1}<0}F_n(xn,s_1,\ldots,s_{n-1},0) ds_1\ldots ds_{n-1} 
=\ga(x)
\end{eqnarray}
almost surely, and it is shown in Lemma \ref{concz} that $\gamma$ is a concave function on $(-\infty,0)$. Therefore, by (\ref{dag}), Theorem \ref{ino2z} and Laplace's method, we would expect
$$-\Psi(m)=\sup_{x<0}[mx+\gamma(x)]=(-\gamma)^*(m)$$
and hence by inversion, $\gamma=-(-\Psi)^*$.
The free energy density for our model of a directed polymer in a Brownian environment in
$1+1$ dimensions, defined 
in (\ref{fre}), can then be expressed in terms of the digamma function by first using the Brownian scaling property:
\begin{eqnarray*}
	\lim_{n \to \infty} \frac{1}{n} \log Z_n(\beta)&=&\gamma(-\beta^2)-2\log |\beta|\\
	&=&-(-\Psi)^*(-\beta^2)-2\log |\beta|.
\end{eqnarray*}

The heuristic argument above is made rigorous by 
the following Theorem, which is the main result of this paper.

\newtheorem{ino4}[ino]{Theorem}
\begin{ino4} \label{ino4z}Almost surely,
\begin{eqnarray*}\lim_{n \to \infty} \frac{1}{n} \log Z_n(\beta)=
f(\beta)
\end{eqnarray*} 
where $f$ is defined by \eqref{putz},
$\Psi(m) \equiv \Gamma'(m)/\Gamma(m)$ is the restriction of the digamma function to $(0,\infty)$, and
$(-\Psi)^*$ is the convex dual of the function $-\Psi$.
\end{ino4}

\section{Proofs}
\label{flint}

We begin by defining for $t\geq0$ and $x<0$
\begin{eqnarray}
	L_n(t)&=&\sup_{0 \leq s_1 \leq \ldots \leq s_{n-1}\leq t}E_n(0,s_1,\ldots,s_{n-1},t)
	\label{le6z} \\
	\nonumber k_{m,n}(x)&=&\log \int_{\substack{xn<s_1<\ldots\\ \ldots<s_{n-m}<xm}} 
	F_{n+1}^{m+1}(xn,s_1,\ldots,s_{n-m},xm)ds_1\ldots ds_{n-m}\\
\mathcal{Z}_n(x)&=&\int_{xn<s_1<\ldots<s_{n-1}<0}
	F_n(xn,s_1,\ldots,s_{n-1},0)
	ds_1\ldots ds_{n-1}\label{gen}\\ 
\gamma_n(x)&=&\frac{1}{n} \log \mathcal{Z}_n(x).
\end{eqnarray}
and recording the following lemma:
\newtheorem{rec}[ino]{Lemma}
\begin{rec} \label{km}
\begin{enumerate}
\item There exists a constant $c \in \mathbb{R}$ such that 
	$\frac{1}{n}L_n(nt) \to c\sqrt{t}$ almost surely and in expectation.
\item There exists a function $\ga:(-\infty,0) \to \rr$ such that given $x<0$, 
\begin{eqnarray}
	\label{pri} k_{0,n}(x)/n \to \gamma(x) \label{damoz}
\end{eqnarray}
almost surely.
\item The function $\gamma$ is continuous on $(-\infty,0)$.
\item Given $x<0$, $\lim_{n \to \infty}\gamma_n(x)=\gamma(x)$ almost surely.

\end{enumerate}
\end{rec}
\textbf{Proof.}
The first part of the lemma follows from Brownian scaling and Kingman's subadditive ergodic theorem (see \cite{MR1935124} for details). 
For the second part, observe that 
\begin{eqnarray*}
k_{0,m}(x)+k_{m,n}(x)&=&\log \int_{\substack{xn<s_1<\ldots<s_{n}<0\\ s_{n-m}<xm<s_{n-m+1}}} F_{n+1}(xn,s_1,\ldots,s_n,0)ds_1\ldots ds_{n} \\
&\leq& k_{0,n}(x)
\end{eqnarray*}
and $k$ is therefore superadditive for fixed $x$. By construction we have the  required conditions for Kingman's subadditive ergodic theorem, 
and so we may define a function $\gamma:(-\infty,0) \to \rr$ by \eqref{damoz}. For the third part we have by Brownian scaling, 
for $x<0$ and $\delta > x$
\begin{eqnarray}
	k_{0,n}(x-\delta)&=_d& \log \int_{\substack{0<s_1<\ldots \\ \ldots<s_{n}<n}}\exp \big(\sqrt{-x+\delta}E_{n+1}(0,s_1,\ldots,s_n,n)\big) ds_1 \ldots ds_{n}
	 \nonumber\\&& \quad+n \log (-x+\delta)
		\nonumber\\
	&\leq& \log \int_{\substack{0<s_1<\ldots \\ \ldots<s_{n}<n}}\exp \big(\sqrt{-x}E_{n+1}(0,s_1,\ldots,s_n,n)\big) ds_1 \ldots ds_{n}
	 \nonumber\\&& \quad+(\sqrt{-x+\delta}-\sqrt{-x})L_{n+1}(n)+n \log (-x+\delta)
	 \nonumber\\&=_d&
		k_{0,n}(x) + (\sqrt{-x+\delta}-\sqrt{-x})L_{n+1}(n)	
	\nonumber\\ && \quad -n \log (-x) +n \log (-x + \delta)
	\label{starz}
\end{eqnarray}
where $=_d$ denotes equality in distribution. Let $c$ be the constant in the first part of this lemma; then by Brownian scaling and Slutsky's Lemma (see for example \cite{MR1609153}), $\frac{1}{n}L_{n+1}(n) \Rightarrow c$.  
By the second part of this lemma and \eqref{starz}, 
\begin{eqnarray*}
	\gamma(x-\delta)&\leq& \gamma(x)+c (\sqrt{-x+\delta}-\sqrt{-x}) -\log(-x)+
	\log(
	-x+\delta)
\end{eqnarray*}
By symmetry, 
\begin{eqnarray}
V_n:=\inf_{0 \leq s_1 \leq \ldots \leq s_{n-1}\leq n}E_n(0,s_1,\ldots,s_{n-1},n)
=_d -L_n(n) \label{efen}
\end{eqnarray}
and so we have similarly
\begin{eqnarray*}
\gamma(x-\delta)
	&\geq& \gamma(x)-c (\sqrt{-x+\delta}-\sqrt{-x}) -\log(-x)+\log(
	-x+\delta).
\end{eqnarray*}
So
\begin{eqnarray}
	|\gamma(x-\delta)-\gamma(x)|
	&\leq& c |\sqrt{-x+\delta}-\sqrt{-x}| \nonumber\\
	&& \quad +|\log(
	-x+\delta)-\log(-x)|. \label{hallz}
\end{eqnarray}
For the fourth part we have
\begin{eqnarray*}
	k_{0,n}(x)&\leq&\log \left( e^{2\overline{B}(n,0,x)}\int_{xn<s_1<\ldots<s_{n}<0} 
	F_n(xn,s_2,\ldots,s_n,0)
	ds_1\ldots ds_{n} \right) \\
&\leq& \log \left( xne^{2\overline{B}(n,0,x)}\int_{xn<s_2<\ldots<s_{n}<0} F_n(xn,s_2,\ldots,s_n,0) ds_2\ldots ds_{n}\right)
\\
&=&\log \left( xne^{2\overline{B}(n,0,x)}\mathcal{Z}_n(x) \right)
\end{eqnarray*}
where for $y\leq x\leq 0$,
\begin{eqnarray} \label{overl}
	\overline{B}(n,x,y)=\max_{i=n-1,n,n+1}\sup_{yn \leq r<s\leq xn}|B^{(i)}(r,s)|.
\end{eqnarray}
From Borell's inequality and the Borel-Cantelli Lemma, there exists a null set $\mathcal{N}$ such that on its complement $\mathcal{N}^c$, 
$\overline{B}(n,x,y)/n \to 0$ for all $x,y$, so
$$\gamma(x) \leq \liminf_{n \to \infty}\gamma_n(x) \quad \mbox{ a.s. }$$
Now for $\epsilon>0$,
\begin{eqnarray*}
k_{0,n}(x+\epsilon)&\geq&\log \left(
e^{-2\overline{B}(n,0,x+\epsilon)}\int_{\substack{(x+\epsilon)n<s_1<xn \\
xn<s_2<\ldots \\ \ldots<s_{n}<0}}
F_n(xn,s_2,\ldots,s_n,0)
ds_1\ldots ds_{n}\right)
\\
&=&\log \left( \epsilon n e^{-2\overline{B}(n,0,x+\epsilon)}\mathcal{Z}_n(x)\right)
\end{eqnarray*}
so
$$\gamma(x+\epsilon) \geq \limsup_{n \to \infty}\gamma_n(x)\quad \mbox{ a.s. }$$
and the result follows by the third part of this lemma.
\hfill $\Box$

Let $\mathbb{Q}_+=\mathbb{Q} \cap (0,\infty)$, $\mathbb{Q}_-=-\mathbb{Q_+}.$
\newtheorem{inc}[ino]{Lemma}
\begin{inc} \label{incz}
There exists a null set $\mathcal{M}$ such that 
the following statement holds on its complement $\mathcal{M}^c$:

 $\lim_{n \to \infty}\gamma_n(x)=\gamma(x)$
and 
$
\liminf_{n \to \infty}\gamma_n(y) \geq \gamma(x)
$ for every $x \in \mathbb{Q}_-$ and $y<x$.

\end{inc}
\textbf{Proof.} Let $x \in \mathbb{Q}_-$ and $y<x$. From (\ref{gen}) we have
\begin{eqnarray} 
	\mathcal{Z}_n(y) &\geq&\int_{xn<s_1<\ldots<s_{n-1}<0}F_n(yn,s_1,\ldots,s_{n-1},0)ds_1\ldots ds_{n-1}
	\nonumber\\
	&=& e^{B^{(n)}_{(yn,xn)}}\mathcal{Z}_n(x).\label{pork}
\end{eqnarray}
By 
Lemma \ref{km}
 there exists a null set $\mathcal{N}_x$ such that on $\mathcal{N}_x^c$, 
\begin{eqnarray*}
	\lim_{n \to \infty}\gamma_n(x)=\gamma(x).
\end{eqnarray*}
With $\mathcal{N}$ as in the proof of Lemma \ref{km}, let $\mathcal{M}=\mathcal{N} \cup \bigcup_{x \in \mathbb{Q}_-}\mathcal{N}_x$.
\hfill $\Box$

\newtheorem{alm}[ino]{Lemma}
\begin{alm} \label{almz}
Almost surely, $\lim_{n \to \infty}\gamma_n(x)=\gamma(x)$ for all $x<0$.
\end{alm}
\textbf{Proof.}
Choose $x,y \in \mathbb{Q}_-$ with $y<x$. Then if $z \in [y,x]$,
\begin{eqnarray*}
	\mathcal{Z}_n(z)&=&\int_{zn<s_1<\ldots<s_{n-1}<0}
	F_n(zn,s_1,\ldots,s_{n-1},0)
	ds_1\ldots ds_{n-1} 
	\\
	&=&I_1(n,z)+I_2(n,z)
\end{eqnarray*}
where
\begin{eqnarray*}
	I_1(n,z)&=&e^{B^{(n)}_{(zn,xn)}} \mathcal{Z}_n(x)\\
	I_2(n,z)&=&\int_{s_1=zn}^{xn}\int_{\substack{s_1<\ldots \\ \ldots<s_{n-1}<0}}
	F_{n-1}(s_1,\ldots,s_{n-1},0)
	ds_2\ldots ds_{n-1} e^{B^{(n)}_{(zn,s_1)}}ds_1\\
	&& \quad  
	\end{eqnarray*}
Now
$
	I_1(n,z) \leq e^{\overline{B}(n,x,y)}\mathcal{Z}_n(x)
$
and
\begin{eqnarray*}
	I_2(n,z) &\leq&
	(x-z) ne^{2\overline{B}(n,x,y)}\mathcal{Z}_{n-1}\Big(\frac{ny}{n-1}\Big)
\end{eqnarray*}
therefore
\begin{eqnarray} \label{unif}
	\sup_{z \in [y,x]}\mathcal{Z}_n(z) \leq e^{\overline{B}(n,x,y)}\mathcal{Z}_n(x) +
	(x-y)ne^{2\overline{B}(n,x,y)}\mathcal{Z}_{n-1}\Big(\frac{ny}{n-1}\Big).
\end{eqnarray}
 Take $\mathcal{M}$ as in the proof of Lemma \ref{incz}.
Now 
\begin{eqnarray*}
	\mathcal{Z}_{n-1}\Big(\frac{ny}{n-1}\Big)&\geq& e^{-\overline{B}(n,0,y)}\mathcal{Z}_{n-1}(y) 
\end{eqnarray*}
so $
\liminf_{n \to \infty} \frac{1}{n} \log \mathcal{Z}_{n-1}\Big(\frac{ny}{n-1}\Big) \geq \gamma(y)$
on $\mathcal{M}^c$; and if $\epsilon \in \mathbb{Q}_-$ then
\begin{eqnarray*}
	\mathcal{Z}_n(y+\epsilon) &\geq& \int_{\substack{(y+\epsilon)n<s_1<yn \\ yn<s_2<\ldots<s_{n-1}<0}} F_{n-1}(yn,s_2,\ldots,s_{n-1},0)ds_1\ldots ds_{n-1}
	\\&& \quad \times e^{-2\overline{B}(n,0,y+\epsilon)}\\
&=&(-\epsilon) n e^{-2\overline{B}(n,0,y+\epsilon)}\mathcal{Z}_{n-1}\Big(\frac{ny}{n-1}\Big)
\end{eqnarray*}
hence
\begin{eqnarray*}
	\gamma(y+\epsilon) \geq \limsup_{n \to \infty} \frac{1}{n}\log \mathcal{Z}_{n-1}\Big(\frac{ny}{n-1}\Big)
\end{eqnarray*}
on $\mathcal{M}^c,$
and letting $\epsilon \to 0 $ in $\qq$ gives $	\lim_{n \to \infty} \frac{1}{n}\log \mathcal{Z}_{n-1}\Big(\frac{ny}{n-1}\Big)=\gamma(y)$ on $\mathcal{M}^c$.
 Then 
by (\ref{unif}) and the proof of Lemma \ref{incz}, on $\mathcal{M}^c$
\begin{eqnarray}
	\gamma(x) &\leq& \liminf_{n \to \infty} \frac{1}{n}\log 
	\inf_{z \in [y,x]}\mathcal{Z}_n(z) \leq
	\limsup_{n \to \infty} \frac{1}{n}\log \sup_{z \in [y,x]}
	\mathcal{Z}_n(z) \nonumber\\ &\leq& \limsup_{n \to \infty} \frac{1}{n}\log 
	e^{\overline{B}(n,x,y)}\mathcal{Z}_n(x) \vee \nonumber\\ && \quad \limsup_{n \to \infty} \frac{1}{n}\log
	(x-y)ne^{2\overline{B}(n,x,y)}\mathcal{Z}_{n-1}\Big(\frac{ny}{n-1}\Big)
		\nonumber\\
	&=& \gamma(x) \vee \gamma(y)
	=\gamma(y). \label{gon2}
\end{eqnarray}
The result now follows from \eqref{hallz}.
\hfill $\Box$
\newtheorem{conc}[ino]{Lemma}
\begin{conc} \label{concz}
The function $\gamma$ is concave.
\end{conc}
\textbf{Proof.} For $x,y<0$ and $\alpha \in (0,1)$,
\begin{eqnarray}
\ga_n(\al y + (1-\al)x)\geq \frac{[\alpha n]}{n}G_n + \frac{k_n}{n} H_n \label{goo}
\end{eqnarray}	
where
\begin{eqnarray*}
	G_n&=&\frac{1}{[\alpha n]}\log \int_{\substack{(\alpha y+(1-\alpha)x)n<s_1<\ldots \\ \ldots<s_{[\alpha n]-1}<(1-\alpha)x n}} 
	F_{n}^{k_n}((\al y+(1-\al)x)n,s_1,\ldots \\ &&  \qquad \qquad \qquad \ldots,s_{[\al n]-1},(1-\al)xn)
	ds_1\ldots ds_{[\alpha n]-1}
	\\
	k_n&=&n-[\alpha n]+1\\
	H_n&=& \frac{1}{k_n}\log \int_{\substack{(1-\alpha) xn<s_{[\alpha n]}<\ldots \\ \ldots<s_{n-1}<0}}
	F_{k_n}((1-\al)xn,s_{[\al n]},\ldots \\ &&
	\qquad \qquad \qquad \ldots,s_{n-1},0)
	ds_{[\alpha n]}\ldots ds_{n-1}
\end{eqnarray*}
Now
\begin{eqnarray*}
	G_n &=_d&
	\ga_{[\al n]}\bigg(\frac{\alpha n}{[\alpha n]} y\bigg)\\
	H_n &=&
	\ga_{k_n}\bigg(\frac{(1-\alpha) n}{k_n} x\bigg).
\end{eqnarray*}
Choose $w,u\in \mathbb{Q}_-$ with $w<x<u$, and choose $\epsilon, \delta>0$. Then $\exists \; n_0 \in \mathbb{N}$ such that $\forall \; n \geq n_0$, $w<\frac{(1-\alpha)n}{k_n}x<u$. Also, by (\ref{gon2}) $\exists \; n_1 \in \mathbb{N}$ such that $\forall \; n \geq n_1$, 
\begin{eqnarray*}
	P\Big(\gamma(u)-\epsilon \leq \frac{1}{n}\log 
	\inf_{z \in [w,u]}\mathcal{Z}_n(z) \leq \frac{1}{n}\log \sup_{z \in [w,u]}
	\mathcal{Z}_n(z) 
	\leq \gamma(w)+\epsilon\Big)>1-\delta
\end{eqnarray*}
Therefore if $k_n\geq n_0 \vee n_1$,
\begin{eqnarray*}
	P(\gamma(u)-\epsilon\leq H_n \leq \gamma(w)+\epsilon)>1-\delta
\end{eqnarray*}
and therefore $H_n \Rightarrow \gamma(x)$. Similarly $G_n \Rightarrow \gamma(y)$, and hence by Slutsky's Theorem 
\begin{eqnarray*}
	\frac{[\alpha n]}{n}G_n + \frac{k_n}{n} H_n \Rightarrow \alpha \gamma(y)+(1-\alpha) \gamma(x)
\end{eqnarray*}
Hence by (\ref{goo}),
\begin{eqnarray*}
		\gamma(\alpha y+(1-\alpha)x)&\geq&\alpha \gamma(y)+(1-\alpha)\gamma(x).
\end{eqnarray*}
\hfill $\Box$

\textbf{Proof of Theorem \ref{ino4z}}
Choose $m>0$. Define a probability density function on $(-\infty,0)$ by
$$\kk_n^{(m)}(dx)=\frac{1}{\Xi_n(m)}\exp n(mx+\gamma_n(x))dx
$$
where
$$\Xi_n(m)=\int_{-\infty}^0 \exp n(mx+\gamma_n(x))dx.
$$
In the nomenclature of statistical physics, 
$\kk_n^{(m)}$ is the \emph{Kac density} and $\Xi_n(m)$
is the \emph{grandcanonical partition function}.
Theorem \ref{ino2z} says that given $\theta > -m,$
the convergences
\begin{eqnarray}
	\frac{1}{n}\log \Xi_n(m) &\to & -\Psi(m) \label{gripz}\\
\frac{1}{n}\log \int_{-\infty}^0e^{n \theta x}\kk_n^{(m)}(dx)
&\to & \Psi(m)-\Psi(m+\theta)=:\Lambda_m(\theta)\nonumber \end{eqnarray}
hold almost surely as $n \to \infty$. 
Choose $\ep>0$; then since $\Lambda_m^*(\Lambda'_m(0)+\ep)$ and $\Lambda_m^*(\Lambda'_m(0)-\ep)$ are strictly positive, we can apply the Chernoff bound  to give
$\kk_n^{(m)}((\Lambda'_m(0)-\ep,\Lambda'_m(0)+\ep)) \to 1$ almost surely as $n \to \infty$. Letting $\ep \to 0$ in $\qq$ 
gives that $\kk_n^{(m)}$ is almost surely concentrated on $\Lambda_m'(0)=-\Psi'(m)$ as $n \to \infty$. Therefore for any $x \in \mathbb{Q}_-$ and $\ep>0$, using \eqref{pork} we have 
\begin{eqnarray*}
	\Xi_n(m) &\geq& \int_{x-\ep}^x \exp n(my+\ga_n(y))dy\\
	&\geq& \exp \{n(m(x-\ep)+\ga_n(x))\}\int_{x-\ep}^x \exp 
	B^{(n)}_{(yn,xn)}dy\\
	&\geq&\ep \exp \{n(m(x-\ep)+\ga_n(x))-\overline{B}(n,x,x-\ep)\}
\end{eqnarray*}
Therefore by Lemma \ref{almz} and \eqref{gripz},
\begin{eqnarray*}
	-\Psi(m) &\geq& m(x-\ep)+\ga(x)
\end{eqnarray*}
and we may let $\ep \to 0$ and appeal to the regularity of $\Psi$ to conclude that  for all $x<0$
$$\ga(x) \leq \inf_{m \in \qq_+}(-mx-\Psi(m))=-(-\Psi)^*(x).$$
For the reverse inequality we note that for $x \in \qq_-$ and 
$\ep \in (0,-x)$
 we may
choose $m>0$ such that $-\Psi'(m)=x+\ep$ (see for example \cite{MR1225604,MR2039096}). Then by \eqref{pork},
\begin{eqnarray*}
	\kk_n^{(m)}((x,x+2\ep)) &=& \frac{1}{\Xi_n(m)}\int_x^{x+2\ep} \exp
	\{n(my+\ga_n(y))\}dy\\
	&\leq & \frac{2\ep}{\Xi_n(m)} \exp\{
	n(m(x+2\ep)+\ga_n(x))-\overline{B}(n,x+2\ep,x)\}
\end{eqnarray*}
Therefore using Lemma \ref{almz},
\begin{eqnarray*}
	0&=& \lim_{n \to \infty}\frac{1}{n} \log \kk_n^{(m)}((x,x+2\ep))
	\leq m(x+2\ep)+\ga(x)+\Psi(m)
\end{eqnarray*}
almost surely, in which case we may let $\ep \to 0$ 
to conclude that 
$$\ga(x) \geq \inf_{m >0}(-mx-\Psi(m))=-(-\Psi)^*(x).
$$
Therefore $\ga(x)=-(-\Psi)^*(x)$. Finally \eqref{zed},
the last part of Lemma \ref{km} and Brownian scaling allow us to conclude that almost surely,
\begin{eqnarray*}\lim_{n \to \infty}\frac{1}{n} \log Z_n(\beta)=\left\{ \begin{array}{ll}
	\gamma(-\beta^2)-2\log |\beta| &: \beta \neq 0 \\
	1&: \beta =0 \end{array}\right.
\end{eqnarray*}
as required. \hfill $\Box$


\section{Analyticity of the free energy density and a large deviation principle}
\label{an}
\newtheorem{ana}[ino]{Theorem} 
\begin{ana} \label{anaz}
The function $f$ defined in \eqref{putz}
 is analytic and strictly convex on $\mathbb{R}$, $f'(0)=0$, and $\lim_{\beta\to \infty}f(\beta)/\beta=2.$
\end{ana}
\textbf{Proof} For $x<0$ we have
\begin{eqnarray*}
(-\Psi)^*(x)=\sup_{\theta>0}[x\theta+\Psi(\theta)]
\end{eqnarray*}
Denoting by $\Psi_n$ the $n$th derivative of the function $\Psi$, and noting that $\Psi_2$ is strictly negative everywhere (see for example \cite{MR2039096}), we have
\begin{eqnarray*}
	(-\Psi)^*(x)=x\Psi_1^{-1}(-x)+\Psi(\Psi_1^{-1}(-x))
\end{eqnarray*}
and since $\Psi_1$ is an invertible analytic function with nonzero derivative, its inverse is analytic. 
Therefore 
$f$ is analytic everywhere except possibly at 0.
 
To investigate the behaviour of $f$ near $0$, let $a=\Psi_1^{-1}(\beta^2)$. Then
\begin{eqnarray} \label{then}
	f(\beta)=a\Psi_1(a)-\Psi(a)-\log \Psi_1(a)
\end{eqnarray}
Now $a \to \infty$ as $\beta \to 0$, and from \cite{MR1225604} we have
\begin{eqnarray} \label{dont}
	\Psi_1(x) \sim \frac{1}{x} + \frac{1}{2x^2} + \frac{1}{6x^3} \quad (x \to \infty)
\end{eqnarray}
therefore recalling (\ref{coff}),
\begin{eqnarray*}
	f(\beta) = 1 + O(a^{-1}) \quad (\beta \to 0)
\end{eqnarray*}
One further application of (\ref{dont}) to this expression gives $f'(0)=0$.

Since $f-1$ is an asymptotic logarithmic moment generating function (see the proof of Lemma \ref{anoldpz} below), $f$ is convex. 
The regularity of $\Psi$ implies further that $f$ is strictly convex.

The power series
\begin{eqnarray*}
	\log \Gamma(1+z)=-\log(1+z)+z(1-\xi)+\sum_{n=2}^{\infty}(-1)^n[\zeta(n)-1]z^n/n
\end{eqnarray*}
where $\xi$ is Euler's constant and $\zeta$ is the Riemann Zeta Function, is valid for $|z|<2$ \cite{MR1225604}; therefore, since $\Psi(z)=\frac{d}{dz}\log \Gamma(z)$,
it may be differentiated to give power series for $\Psi$ and $\Psi_1$ in a neighbourhood of the origin. Substituting in (\ref{then}) and letting 
$a \to 0$, we obtain $\lim_{\beta\to \infty}f(\beta)/\beta=2.$
\hfill $\Box$

\newtheorem{anoldp}[ino]{Lemma}
\begin{anoldp} \label{anoldpz}
Let $\tau_1 \leq \ldots \leq \tau_{n-1}$ 
be the order statistics for $n-1$ independent random variables having the uniform distribution on the interval $[0,n]$.
Almost surely, conditional on $\mathcal{B}$, the random variable $\frac{1}{n}E_n
(0,\tau_1,\ldots,\tau_{n-1},n)$ satisfies a large deviation principle with rate function $(f-1)^*$.
\end{anoldp}
\textbf{Proof.} 
Choose $\be \in \rr$; then by \eqref{zed}
\begin{eqnarray*}
	E[\exp(\beta E_n
(0,\tau_1,\ldots,\tau_{n-1},n))|\mathcal{B}]&=&\frac{(n-1)!}{n^{n-1}}Z_n(\be)
	\end{eqnarray*}
hence by Stirling's formula and Theorem \ref{ino4z}, 
\begin{eqnarray}
	\lim_{n \to \infty} \fn{1} \log E[\exp(\beta E_n
(0,\tau_1,\ldots,\tau_{n-1},n))|\mathcal{B}]=f(\be)-1 \label{libby}
\end{eqnarray}
 almost surely.
Now if $\al<\nu$ then
$	Z_n(\al) \leq e^{-(\nu-\al) V_n} Z_n(\nu)\label{use}$, where $V_n$ was defined in \eqref{efen};
hence for $\al,\be \in \qq$ with $\beta > \alpha$
\begin{eqnarray*}
	f(\al) &\leq& \liminf_{n \to \infty} \fn{1} \log \inf_{\nu \in (\al,\be)}Z_n(\nu)
	+(\be-\al)c \nonumber\\
	&\leq& \limsup_{n \to \infty} \fn{1} \log \sup_{\nu \in (\al,\be)}Z_n(\nu)
	+(\beta-\alpha)c \leq f(\be) + 2(\be-\al)c
\end{eqnarray*}
almost surely, where $c$ was defined in Lemma \ref{km}. Therefore by the continuity of $f$, there exists a null set $\mathcal{N}$ such that on $\mathcal{N}^c$,
the convergence in \eqref{libby} holds for all $\be \in \rr$. 
\hfill $\Box$

\section{Connection with random matrices and a Brownian directed percolation problem}
\label{appl}

It was shown in \cite{MR1935124} that 
with $L_n$ defined as in \eqref{le6z},
for each $t \geq 0$,
\begin{eqnarray*}
	\lim_{n \to \infty} \frac{1}{n}L_n(nt) = 2\sqrt{t}
\end{eqnarray*}
almost surely. In the notation of Lemma \ref{km} this says that $c=2$.
This result can also be deduced from
random matrix theory using the fact \cite{MR1818248, MR1830441} that $L_n(1)$
has the same law as the largest eigenvalue of a $n \times n$
GUE random matrix. We note here that the present work allows us to deduce the inequality $c \geq 2$.  Since $L_n(n) \geq E_n(0,t_1,\ldots,t_{n-1},n)$
we have that given $\beta$,
\begin{eqnarray*}
	\frac{1}{n}L_n(n) &\geq& \frac{1}{\beta}\frac{1}{n}\log E[\exp(\beta E_n
(0,\tau_1,\ldots,\tau_{n-1},n))|\mathcal{B}]
\end{eqnarray*}
almost surely. Letting $n \to \infty$ and then $\beta \to \infty$, using \eqref{libby} and Theorem \ref{anaz} gives the required inequality.

\bibliography{free2}

\begin{thebibliography}{10}

\bibitem{MR1225604}
Milton Abramowitz and Irene~A. Stegun, editors.
\newblock {\em Handbook of mathematical functions with formulas, graphs, and
  mathematical tables}.
\newblock Dover Publications Inc., New York, 1992.
\newblock Reprint of the 1972 edition.

\bibitem{MR2039096}
Horst Alzer.
\newblock Sharp inequalities for the digamma and polygamma functions.
\newblock {\em Forum Math.}, 16(2):181--221, 2004.

\bibitem{MR1818248}
Yu. Baryshnikov.
\newblock G{UE}s and queues.
\newblock {\em Probab. Theory Related Fields}, 119(2):256--274, 2001.

\bibitem{MR1939654}
Philippe Carmona and Yueyun Hu.
\newblock On the partition function of a directed polymer in a {G}aussian
  random environment.
\newblock {\em Probab. Theory Related Fields}, 124(3):431--457, 2002.

\bibitem{MR2117626}
Francis Comets and Nobuo Yoshida.
\newblock Brownian directed polymers in random environment.
\newblock {\em Comm. Math. Phys.}, 254(2):257--287, 2005.

\bibitem{MR1043640}
B.~Derrida.
\newblock Directed polymers in a random medium.
\newblock {\em Phys. A}, 163(1):71--84, 1990.
\newblock Statistical physics (Rio de Janeiro, 1989).

\bibitem{MR1833627}
Daniel Dufresne.
\newblock An affine property of the reciprocal {A}sian option process.
\newblock {\em Osaka J. Math.}, 38(2):379--381, 2001.

\bibitem{MR1825324}
Daniel Dufresne.
\newblock The integral of geometric {B}rownian motion.
\newblock {\em Adv. in Appl. Probab.}, 33(1):223--241, 2001.

\bibitem{MR1609153}
Richard Durrett.
\newblock {\em Probability: theory and examples}.
\newblock Duxbury Press, Belmont, CA, second edition, 1996.

\bibitem{MR1830441}
Janko Gravner, Craig~A. Tracy, and Harold Widom.
\newblock Limit theorems for height fluctuations in a class of discrete space
  and time growth models.
\newblock {\em J. Statist. Phys.}, 102(5-6):1085--1132, 2001.

\bibitem{MR1935124}
B.~M. Hambly, James~B. Martin, and Neil O'Connell.
\newblock Concentration results for a {B}rownian directed percolation problem.
\newblock {\em Stochastic Process. Appl.}, 102(2):207--220, 2002.

\bibitem{MR1696208}
Hiroyuki Matsumoto and Marc Yor.
\newblock A version of {P}itman's {$2M-X$} theorem for geometric {B}rownian
  motions.
\newblock {\em C. R. Acad. Sci. Paris S\'er. I Math.}, 328(11):1067--1074,
  1999.

\bibitem{MR1833628}
Hiroyuki Matsumoto and Marc Yor.
\newblock A relationship between {B}rownian motions with opposite drifts via
  certain enlargements of the {B}rownian filtration.
\newblock {\em Osaka J. Math.}, 38(2):383--398, 2001.

\bibitem{MR1865759}
Neil O'Connell and Marc Yor.
\newblock Brownian analogues of {B}urke's theorem.
\newblock {\em Stochastic Process. Appl.}, 96(2):285--304, 2001.

\end{thebibliography}
\end{document}